\documentclass[draft]{article} 
\usepackage[cp1251]{inputenc}         
\usepackage{amssymb,amsfonts,amsmath} 
\usepackage{amsthm}
\newtheorem{proposition}{Proposition}
\newtheorem{theorem}{Theorem}
\newtheorem{lemma}{Lemma}
\newtheorem{corollary}{Corollary}
\newtheorem{remark}{Remark}
\newtheorem{example}{Example}

\newcommand{\cP}{{\mathcal P}}

\newcommand{\bbM}{{\mathord{\mathbb M}}}
\newcommand{\bbR}{{\mathord{\mathbb R}}}

\newcommand{\bbZ}{{\mathord{\mathbb Z}}}
\newcommand{\one}{\mathord{\mathbf1}}

\newcommand{\scal}[2]{{\left\langle{#1},{#2}\right\rangle}}

\newcommand{\Aff}{\mathop{\rm Aff}\nolimits}

\newcommand{\OO}{\mathop{\rm O}\nolimits}

\newcommand{\GL}{\mathop{\rm GL}\nolimits}

\newcommand{\rank}{\mathop{\mathrm{rank}}\nolimits}

\newcommand{\spp}{\mathop{\mathrm {sp}}\nolimits}

\newcommand{\be}{\beta}
\newcommand{\ga}{\gamma}
\newcommand{\Ga}{\Gamma}
\newcommand{\al}{\alpha}

\newcommand{\la}{\lambda}
\newcommand{\La}{\Lambda}

\newcommand{\ka}{\kappa}
\newcommand{\ze}{\zeta}

\newcommand{\wh}[1]{\widehat{#1}}

\let\td=\tilde
\title{On Geometry of Flat Complete Strictly Causal Lorentzian Manifolds}
\author{V.M. Gichev, E.A. Meshcheryakov}
\date{}
\begin{document}
\maketitle
\begin{abstract}
A flat complete causal Lorentzian manifold is called {\it strictly
causal} if the past and the future of each its point are closed
near this point. We consider strictly causal manifolds with
unipotent holonomy groups and assign to a manifold of this type
four nonnegative integers (a signature) and  a parabola in the
cone of positive definite matrices. Two manifolds are equivalent
if and only if their signatures coincides and the corresponding
parabolas are equal (up to a suitable automorphism of the cone and
an affine change of variable). Also, we give necessary and
sufficient conditions, which distinguish parabolas of this type
among all parabolas in the cone.
\end{abstract}

\section{Introduction}
A flat complete Lorentzian manifold can be realized as a quotient
space $M=\bbM_n/\Ga$, where $\Ga\cong\pi_1(M)$ is a discrete
subgroup of Poincare group $\cP_n$ of affine automorphisms of
Minkowski spacetime $\bbM_n$ which acts on $\bbM_n$ freely and
properly. There are several equivalent definitions: $M$ is a
geodesically complete Lorentzian manifold with vanishing torsion
and curvature; $M$ admits an atlas of coordinate charts in
$\bbM_n$ with coordinate transformations in $\cP_n$ such that any
affine mapping of a segment in $\bbR$ to $M$ extends to an affine
mapping of $\bbR$ (a complete affine manifold with a compatible
Lorentzian metric). The Lorentzian metric defines a pair of closed
convex round cones in each tangent space $T_pM$. Choosing one of
them, we get locally a cone field. It can be extended to a global
cone field on $M$ or on a two-sheet covering space of $M$; we
assume that the field is defined on $M$. This is equivalent to the
assumption that linear parts of transformations in $\Ga$ do not
transpose the cones of the past and the future in $\bbM_n$. If $M$
admits no closed timelike curves, then it is called {\it causal}.
For $\Ga$, this means that the orbit of any point $v$ does not
meet the cone at $v$. Fixing origin $o\in\bbM_n$, we may identify
$\bbM_n$ with a real vector space $V$ endowed with a Lorentzian
metric $\ell$ of the signature $(+,-\dots,-)$; the causal
structure is determined by a cone $C$ (one of the two cones
defined by the inequality $\ell(v,v)\geq0$). We say that $M$ is
{\it strictly causal\/} if the past and the future of any point
$p\in M$ are closed near $p$ (they are not closed globally in
general).  In paper \cite{GM}, these manifolds were found up to
finite coverings. Here is a brief description: $M$ is the total
space of a topologically trivial vector bundle with a bounded
linear holonomy group, whose base is $\bbM_k/\Gamma$, $k\leq n$,
where the action of $\Ga$ on $\bbM_k$ is affine and unipotent.
This reduces the problem to the case of unipotent $\Ga$. There are
two types of them. The first (elliptic) is trivial: $\Ga$ is the
group of translations by vectors in a unform lattice in a linear
subspace $T\subset V$ such that $T\cap C=0$. We consider the
second (parabolic) type. The least dimension of such a manifold is
4; moreover, all these manifolds are mutually homothetic in
dimension 4. In this paper, we assign to a manifold of this type
(of any dimension) a parabola in the cone $P_T$ of positive
definite quadratic forms on a Euclidean space $T$. It is
parametrized by a quadratic polynomial with matrix coefficients.
The manifold $M$ is completely determined by this parabola
(considered up to automorphisms of the cone) and a lattice
$\Ga\subset T$. The parabolas corresponding to the manifolds do
not exhaust all parabolas in the cone; we characterize them and
describe some their invariants.

The object of this paper lies in the intersection of the two
well-explored fields: flat complete affine manifolds (see recent
surveys \cite{Abels} and \cite{CDGM03}) and the causality in the
Lorentzian manifolds; the studying of the latter was mainly
stimulated by general relativity (see [4, 5]). In their common
part, we know only some articles of A. D. Alexandrov’s
chronogeometric school (see \cite{Guts} for references) and papers
\cite{GM}, \cite{Fried}, \cite{Barbot}. The paper \cite{Fried}
contains a characterization of two-ended causal 4-manifolds which
can be realized as $H/\Ga$, where $H$ is a subgroup of Poincare
group whose action on Minkowski spacetime is simply transitive,
$\Ga$ is its discrete subgroup. Most of them are not strictly
causal while the nontrivial strictly causal manifolds are never
globally hyperbolic (\cite{Barbot},\cite{GM}); the latter class of
flat Lorentzian (not necessarily complete) manifolds was
considered in \cite{Barbot}. This paper continues article
\cite{GM}, where some parametric description was found for the
flat complete strictly causal Lorentzian manifolds (we give it
below).

\section{Preliminaries and statement of results}
Fixing the origin $o\in\bbM_n$, we realize $\bbM_n$ as a real
vector space $V$ equipped with a Lorentzian form $\ell$ of the
signature $(+,-,\dots,-)$. The set $\ell(v,v)\geq0$ is the union
of two closed convex round cones in $V$. Let $C$ be one of them.
The group $\Ga$ acts freely and properly in $V$ by affine
transformations whose linear parts preserve $\ell$ and $C$
invariant. We denote by $\ka$ the quotient mapping $\bbM_n\to
M=V/\Ga$ and define the {\it past} $P_p$ and the {\it future}
$F_p$ of $p\in M$ by
\begin{eqnarray*}
P_p=\ka(v-C),\qquad F_p=\ka(v+C),\qquad v\in\ka^{-1}(p).
\end{eqnarray*}
Clearly, $P_p$ and $F_p$ do not depend on the choice of $v$.
Hence, the projection of the cone field to $M$ is well-defined:
\begin{eqnarray*}
&C_p=d_v\ka(v+C)\subset T_pM,\quad v\in \kappa^{-1}(p).
\end{eqnarray*}
The manifold $M$ is said to be {\it causal} if $M$ admits no
closed piecewise smooth timelike paths. A smooth path $\eta$ is
called {\it timelike} if $\eta'(t)\in C_{\eta(t)}$ for all $t$;
for {\it lightlike} paths, $\eta'(t)\in
\partial C_{\eta(t)}$ (note that lightlike paths are timelike).
Clearly, each timelike curve in $M$ can be lifted to a timelike
curve in $V$ and the projection into $M$ of a timelike curve in
$V$ is timelike. We say that an isometry of Lorentzian manifolds
is {\it causal} if it preserved the orientation of timelike
curves. Then the manifolds are said to be {\it causally
isometric.} An {\it affine manifold\/} is a manifold with affine
coordinate transformations. A {\it complete} affine manifold can
be realized as $V/\Ga$ where $\Ga$ is a free and proper group of
affine transformations of $V$; in our setting, they belong to
$\cP_n$. Given $\ga\in\Ga$, put
\begin{eqnarray}
&\ga(v)=\la(\ga)v+\tau(\ga), \quad\mbox{where}\quad
 \la(\ga)\in\OO(\ell),\ \tau(\ga)\in V;\label{affde}\\
&G=\la(\Ga),\nonumber
\end{eqnarray}
where $\OO(\ell)$ denotes the group of all linear transformations
of $V$ preserving $\ell$. Clearly, the mapping
$\la:\,\Ga\to\OO(\ell)$ is a homomorphism and for all
$\ga,\nu\in\Ga$
\begin{eqnarray*}
&\tau(\ga\circ\nu)=\la(\ga)\tau(\nu)+\tau(\ga).
\end{eqnarray*}
According to \cite[Theorem~1]{GM}, each strictly causal flat
complete Lorentzian mani\-fold is finitely covered by  the total
space of a vector bundle with (arbitrary) bounded holonomy group
and unipotent base (the latter means that $\Ga$ consists of affine
transformations with unipotent linear parts). Thus we consider
only the unipotent case. By \cite[Theorem~2]{GM}, a unipotent
manifold of this type can be finitely covered by a manifold
described below.

{\bf The Main Construction}. Let $v_0,v_1\in\partial C$ satisfy
$\ell(v_0,v_1)=1$ and put
\begin{eqnarray}
&L=\bbR v_0,\quad W=L^\bot, \quad N=W\cap v_1^\bot;\label{lwn}\\
&l_0(v)=\ell(v_0,v).\nonumber
\end{eqnarray}
The hyperplane $W=N\oplus L$ is tangent to $\partial C$ at $v_0$,
while the set $W\cap C$ is a ray. The form $\ell$ is nonpositive
and degenerate in $W$ and nondegenerate and negative in $N$. Let
$T\subseteq N$ be a linear subspace and put $\td T=T+L$. We will
often identify $T$ with $\td T/L$. Let $\Ga$ be a lattice
(cocompact discrete subgroup of the additive group) in $T$ and $a$
be an $\ell$-symmetric linear mapping:
\begin{eqnarray}
&a:\,T\to N,\nonumber\\
&\ell(ax,y)=\ell(x,ay), \quad x,y\in T.\label{asymm}
\end{eqnarray}
The affine action $x\to\ga_x$ of $T$ in $V$ is defined by formulas
\begin{eqnarray}
&\la(x)v = v+l_0(v)ax-\left(\ell(ax,v)+\frac12\,l_0(v)\ell(ax,ax)\right)v_0,\label{forla}\\
&\tau(x)=x-\frac12\,\ell(ax,x)v_0;\label{forta}\\
&\ga_x(v)=\la(x)v+\tau(x).\nonumber
\end{eqnarray}
The following condition is necessary and sufficient for the action
of $T$ to be free and for action of  $\Ga$ to be free and proper
\cite[Lemma~19]{GM}:
\begin{equation}\label{freea}
\ker({\mathbf1}+sa)=0\quad\mbox{for all}\quad s\in\bbR,
\end{equation}
where $\mathbf1$ is the identical mapping.\qed

The quotient mappings $V\to V/L$ and $V\to V/\Ga$ are denoted by
$\phi$ and $\ka$, respectively. Some simple properties of the
above action are stated in next lemma.

\begin{lemma}\label{trivi}
If (\ref{freea}) is true then the following hold for the action
(\ref{lwn})--(\ref{forta}).
\begin{itemize}
 \item[{(1)}]
If $x\in T$ and $ax\neq0$, then the line $L$ is precisely the set
of all fixed points of $\la(x)$ in $C\cup(-C)$; translations by
vectors in $L$ commute with $\ga_x$ for all $x\in T$. \item[(2)]
The action of $\Ga$ in the quotient space $V/L$ is free and
proper. Every hyperplane
\begin{eqnarray*}
W_s=\{v\in V:\,l_0(v)=s\}\subset V
\end{eqnarray*}
is $\Ga$-invariant, and $\Ga$ acts by pure translations in
$W_s/L\subset V/L$. \item[(3)] The mapping $\phi$ is one-to-one on
every $T$-orbit in $V$.\qed
\end{itemize}
\end{lemma}
The set of common fixed points of $G=\la(\Ga)$ in $V$ may be
greater than $L$. If $a=0$ then $\Ga$ and $T$ act by translations:
$\ga_x(v)=v+x$. In \cite{GM} this case was called {\it elliptic}
and considered separately. In this paper, we combine elliptic and
{\it parabolic} ($a\neq0$) cases.

The affine structure makes it locally possible to decide whether
two vectors are parallel or not. Hence, the parallel transport of
vectors along curves is well-defined. Applying this to loops at
$p\in M$, we obtain the holonomy representation
$\la_p:\,\pi_1(M,p)\cong\Ga\to\GL(T_pM)$. We have the natural
identification $\la_p(\pi_1(M,p))=\la(\Ga)=G$. The following
theorem is an observation that refines \cite[Theorem~2]{GM}, where
an analogous assertion was proved up to finite coverings and
without mentioning of holonomy. A linear group is said to be {\it
unipotent} if, in some linear base, it can be realized by
triangular matrices whose diagonal entries are equal to 1. We say
that $\Ga\subset\Aff(V)$ is unipotent if $G=\la(\Ga)$ has this
property (it can be considered as a unipotent linear group in the
space $V\oplus\bbR$).
\begin{theorem}\label{uniho}
A flat complete strictly causal Lorentzian manifold admits a
realization above if and only if its holonomy group is unipotent.
\end{theorem}
We say that manifolds of Theorem~\ref{uniho} are {\it unipotent}.
\begin{corollary}
The fundamental group of a unipotent manifold is isomorphic to
$\bbZ^m$.\qed
\end{corollary}
Let $\wh\Ga$ denote the algebraic (i.e. in the Zariski topology)
closure of $\Ga$ in the group $\Aff(V)$ of all affine
transformations of $V$ (clearly, $\wh\Ga\subset\cP(V)$).
\begin{proposition}\label{zarcl}
The algebraic closure $\wh\Ga$ of $\Ga$ coincides with the image
of $T$ under the embedding to $\cP(V)$ defined by (\ref{forla})
and (\ref{forta}); in particular, $\wh\Ga$ is isomorphic to the
vector group $T$ and $\wh\Ga/\Ga$ is a torus. Moreover, all orbits
of $\wh\Ga$ in $V$ are Zariski closed.
\end{proposition}
We assume that $\wh\Ga$ is a subgroup of $\cP_n$ and $T$ is a
linear subspace of $V$ in the sequel. The orbits of $\wh\Ga$ in
$M$ are tori $\wh\Ga/\Ga$. If $a=0$ then they are affine
submanifolds of $M$; otherwise, there is the compatible affine
structure on line bundles over these tori, with lines parallel to
$L$ (see Lemma~\ref{trivi}, (2)). The torus $\wh\Ga/\Ga$ acts
freely on $M$ in both cases, $M/\wh\Ga$ is homeomorphic to a
vector space and $M$ is homeomorphic to the product
$\wh\Ga/\Ga\times M/\wh\Ga$ (\cite[Theorem~2]{GM}).

The manifold $M$ is completely  determined by parameters
$v_0,v_1,T,a,\Ga$. They are not independent (for example, $T$ is
the linear span of $\Gamma$). The vectors $v_0$ and $v_1$ define
$W,N,L$ by (\ref{lwn}). If $v_0$ is fixed then every choice of a
complementary to $L$ subspace $N$ in $W$ determines $v_1$
uniquely: the space $N^\bot$ is two dimensional and intersects
$\partial(C\cup(-C))$ by two lines, $L$ and $\bbR v_1$, and the
position of $v_1$ in the second line is defined by
$\ell(v_0,v_1)=1$. We write
\begin{eqnarray*}
M=M(v_0,v_1,T, a,\Ga)=M(v_0,N,T,a,\Ga)
\end{eqnarray*}
omitting parameters sometimes. Most of them have a natural
geometrical meaning; the following proposi\-tion clarifies it. Let
$M$ be as in Theorem~\ref{uniho}, $p\in M$, $\Ga=\pi_1(M,p)$. We
define some subspaces of $T_pM$ in the notation which agrees with
the main construction:
\begin{itemize}
\item[] $T$: the tangent space to the orbit $\wh\Ga p\,$; \item[]
$H$: the linear span of vectors $(\la_p(x)-\one)v$, where
$x\in\Ga$, $v\in T_pM$; \item[] $U=T+H$, \quad $L=H\cap
H^\bot$,\quad $W=L^\bot$.
\end{itemize}
If $M$ is elliptic then $H=L=0$, $U=T$, $W=V$. Since the action of
$\wh\Ga/\Ga$ in $V$ is free, $T$ may be identified with the Lie
algebra of the torus $\wh\Ga/\Ga$. This defines the exponential
mapping $\exp:\,T\to\wh\Ga p$. For a flat complete affine manifold
$M$ and each $p\in M$, another exponential mapping
$\exp_p:\,T_pM\to M$ is uniquely defined by following conditions:
\begin{itemize}\item[(1)]
$\frac{d}{dt}\exp_p(tu)|_{t=0}=u$ for each $u\in T_pM$,\item[(2)]
 the mapping  $t\to\exp_p(tu)$ is affine.
\end{itemize}
The two exponential mappings are not equal  but their
$\phi$-projections coincide on $T$ by Lemma~\ref{trivi}, (2). Let
$\pi_X$ denote the $\ell$-orthogonal projection to a subspace $X$
(the definition is sound if $\ell$ is nondegenerate in $X$).
\begin{proposition}\label{geome}
Let $M$ be a unipotent nonelliptic manifold, $p\in M$, and $T$,
$H$, $U$, $L$, $W$ be as above. Then $\dim L=1$ and $L$ consists
of fixed points of the holonomy representation.
\begin{itemize}
\item[(1)] Every choice of a generating vector $v_0\in\partial
C_p$ for $L$, an isotropic vector $v_1\perp U$ such that
$\ell(v_0,v_1)=1$, defines the action (\ref{lwn})--(\ref{forta})
of $T$  on the setting $N=v_1^\bot\cap W$ by
\begin{eqnarray*}
ax=\pi_N\la(\exp(x))v_1,\quad x\in T,
\end{eqnarray*}
with $T_pM$ as $V$. \item[(2)] The mapping $\exp_p$ satisfies the
condition $\exp_p(\ga(v))=\exp_p(v)$ for all $v\in T_pM$ and
identifies $M$ with $T_pM/\Ga$. \item[(3)]Put $E=U^\bot\cap N$.
Then $M=E\times M'$, where $M'=\exp_p(E^\bot)$ is a nonelliptic
unipotent submanifold of $M$.
\end{itemize}
\end{proposition}

Each $\ell$-symmetric mappings has an evident structure: if $a$
satisfies (\ref{asymm}) then it admits the unique decomposition
\begin{eqnarray}\label{deca}
a=a'+a'',
\end{eqnarray}
where $a':\,T\to T$ is the self-adjoint transformation of $T$
corresponding to the symmetric bilinear form $\ell(ax,y)$:
\begin{eqnarray*}
&\quad\ell(ax,y)=\ell(a'x,y)=\ell(x,a'y), \quad x,y\in T,
\end{eqnarray*}
and $a''$ is an arbitrary linear mapping
\begin{eqnarray*}
&a'':\,T\to T^\bot\cap N.
\end{eqnarray*}
Put \begin{eqnarray*} R=a''T.
\end{eqnarray*}
The condition (\ref{freea}) can be rewritten as follows:
\begin{equation}\label{freec}
t\in\bbR,\quad a'x=tx\neq0\quad\Longrightarrow\quad a''x\neq0.
\end{equation}
In other words, $a''$ is nondegenerate in all eigenspaces of $a'$
but $\ker a'$ (note that $a'$  has only real eigenvalues and is
semisimple since $a'$ is self-adjoint). The space
\begin{eqnarray*}
\ker a=\ker a'\cap\ker a''\subseteq T
\end{eqnarray*}
acts by pure translations.  It is a trivial summand for the action
of $T$ but this fails in general for $\Gamma$.

There are three natural steps in the construction of $M$:
\begin{enumerate}
\item[(A)] fix $v_0,v_1$, define $L,W,N$ by (\ref{lwn}) and choose
$T\subseteq N$; \item[(B)] pick an $\ell$-symmetric linear
operator $a':\,T\to T$, and, for each its eigenspace $\La_j$, a
linear operator $a_j'':\,\La_j\to T^\bot\cap N$ (which must be
nondegenerate if $\La_j\neq\ker a'$); and set $a''=\sum_ja_j''$,
$a=a'+a''$; \item[(C)] choose a linear basis for $T$ and define
$\Ga$ as the subgroup of the vector group $T$ generated by it.
\end{enumerate}
The first step provides a frame for the second and the third which
are independent of one another. For example, if $T=N$, then $a=0$
and we get an elliptic manifold; if $\dim R=1$, then $a''$ has
rank 1 and can be nondegenerate in all eigenspaces of $a'$ only if
they are one-dimensional. We say that $a'$ and $M(a)$ have the
{\it simple spectrum} if each eigenvalue of $a'$ has multiplicity
$\leq1$.
\begin{remark}\rm
It is not difficult to construct a unipotent manifold with the
simple spectrum and prescribed eigenvalues and eigenvectors of
$a'$ (any orthonormal linear base in $T$). A manifold with the
simple spectrum is determined up to an isometry by $m$ real
numbers (the spectrum of $a'$) and the Gram matrix of $m$ vectors
($a''$-images of eigenvectors) of rank $r\geq1$ which must have
nonzero diagonal elements (this is not a classification since
these parameters do not distinguish some isometric manifolds).
\end{remark}

The causally isometric  manifolds of this type admit realizations
(\ref{lwn})--(\ref{forta}) with identical parameters. We say that
manifolds are {\it almost causally isometric} if they admit
realizations that differ only on the step (C) of the construction
above.

To all $p\in M$ and $x\in\pi_1(M,p)$, there is a realization of
the loop $x$ as a straight line segment. Precisely, this is the
projection into $M$ of the segment with endpoints $\ga_x(v)$ and
$v$, where $\ka(v)=p$ and $\ga_x\in\Ga$ is the affine
transformation corresponding to $x$. Put
\begin{eqnarray}\label{defqv}
&q_v(x)=-\ell(\ga_x(v)-v,\ga_x(v)-v).
\end{eqnarray}
This is a function on the group $\pi_1(M,p)=\Ga$ (the squared
$\ell$-length of the segment mentioned above). If $a=0$ then
$q_v(x)=-\ell(x,x)$ does not depend on $v$ since $\Ga$ acts by
pure translations. In general, $\Ga$ acts by translations in each
hyperplane $W_s/L$ (Lemma~\ref{trivi}, (2)). Since $W\perp L$,
this means that $q_v$ depends only on $s=l_0(v)$. By a
straightforward calculation with (\ref{forla}), (\ref{forta}) we
obtain
\begin{eqnarray}\label{raspi}
q_v(x)=q_s(x)=-\ell((\one+sa)x,(\one+sa)x),\quad\mbox{where}\quad
s=l_0(v).
\end{eqnarray}
Hence $\{q_v:\,v\in V\}$ is one-parameter family of quadratic
forms on $T$. Since $\ell$ is negative definite on $T$, it follows
from (\ref{freea}) that all forms $q_s$ are positive definite on
$\wh\Ga$. Thus, we arrive at a curve in the cone od positive
definition quadratic forms on $\wh \Gamma$ which is a parametrized
by $s=l_0(v)$ . By (\ref{defqv}) and (\ref{raspi}), the change of
origin is equivalent to the shift of the parameter:
\begin{eqnarray}\label{shifs}
s\to s-s_0,\quad s_0=l_0(\td o),
\end{eqnarray}
where $\td o$ is new origin. For all $t>0$, replacing $v_0,v_1$ by
$tv_0,v_1/t$, respectively, we come to the same formulas with
$tl_0, a/t$ instead of $l_0,a$. This corresponds to the change of
variable
\begin{eqnarray}\label{affch}
s\to ts.
\end{eqnarray}
If $t<0$ then the time reverses. Thus, we will consider the curve
$s\to q_s$ up to orientation-preserving affine changes of the
variable $s$.

We will identify $\Gamma$ with $\bbZ^m$ and $\wh\Gamma$ with
$\bbR^m$. More precisely, let $\iota:\,\bbR^m\to T$ be a linear
isomorphism such that $\iota\,\bbZ^m=\Ga$ and let $\scal{\ \,}{\
}$ be the standard inner product in $\bbR^m$. By (\ref{raspi}),
$q_s$ is quadratic in s. Hence, there exist symmetric $m$-matrices
$A,B,C$ such that
\begin{eqnarray}\label{qsabc}
q_s(x)=\scal{(A+2sB+s^2C)z}{z}\quad\mbox{for all}\ z\in\bbR^m,\ \
\mbox{where}\ x=\iota\,z\in T.
\end{eqnarray}
The inequality $S>0$ ($S\geq0$), where $S$ is a matrix, means that
$S$ is symmetric positive definite (respectively, nonnegative). We
denote the set of positive matrices by $P_m$; it is a homogeneous
space of the group $\GL(m,\bbR)$ acting by
\begin{eqnarray*}
&S\to X^\top SX,
\end{eqnarray*}
where $\ ^\top$  stands for the transposition. Moreover, the
involution $S\to S^{-1}$ defines on $P_m$ the structure of a
symmetric space. The condition
\begin{eqnarray}\label{poabc}
Q(s)=A+2sB+s^2C>0\quad\mbox{for all}\ s\in\bbR
\end{eqnarray}
is necessary (but not sufficient) for the matrix valued quadratic
polynomial $Q$ to satisfy (\ref{qsabc}). It implies $A>0$,
$C\geq0$ (note that $B=C=0$ if $a=0$).
\begin{remark}\rm
The affine span of a generic curve of this type is at most
two-dimensional: the linear subspace parallel to it is spanned by
matrices $B$ and $C$. Thus, it is a parabola, which may degenerate
into a ray if $B=0$ and into a point if $B=C=0$ (the case $C=0$,
$B\neq0$ cannot occur due to (\ref{poabc})). Up to an affine
change of the variable $s$, there is only one quadratic
parametrization of a parabola in a plane. Thus, we come to a
geometrical object, a {\it characteristic curve of $M$.}
\end{remark}
A linear change of the variable $z\in\bbR^m$ defined by a real
$m$-matrix $X\in\GL(m,\bbR)$, induces come translation of this
curve in the symmetric space $P_m$:
\begin{eqnarray}\label{linch}
Q(s)\to X^\top Q(s)X.
\end{eqnarray}
We say that $Q(s)$ is a {\it characteristic polynomial of $M$\/}
and denote it by $Q_M$, considering $Q_M$ up to affine changes of
the variable $s$. Put
\begin{eqnarray}\label{nmr}
n=\dim M,\quad m=\dim T,\quad r=\dim R;\quad k=\dim\ker a;
\end{eqnarray}
the 4-tuple $(n,m,r,k)$ will be called a {\it signature of $M$\/}.
It follows from Proposition~\ref{geome} that the signature does
not depend on the realization of $M$ in the form
(\ref{lwn})--(\ref{forta}). Clearly, these numbers satisfy
inequalities
\begin{eqnarray*}\label{ineq}
m+r+2\leq n,\quad r+k\leq m.
\end{eqnarray*}
\begin{theorem}\label{main}
Let manifolds $M$ and $\td M$ be as in Theorem~\ref{uniho}. They
are causally isometric if and only if their signatures coincide
and
\begin{eqnarray}\label{eqcur}
Q_M(s)=X^\top Q_{\td M}(\al s+\be)X
\end{eqnarray}
for some $X\in\GL(m,\bbZ)$, $\al>0$, $\be\in\bbR$.
\end{theorem}
In other words, the manifolds with equal signatures are isometric
if and only if projections of their characteristic curves to
$P_m/\GL(m,\bbZ)$ coincide. Note that $P_m/\GL(m,\bbZ)$ is the
modulus space for Euclidean structures in $m$-tori.  It appears
naturally since Euclidean structures induced by $-\ell$ in the
orbits of the torus $\wh\Ga/\Ga$ in $M/L$ depends only on
$s=l_0(v)$ (the group $L\cong\bbR$ acts on $M=V/\Ga$ since
translations by vectors in $L$ commute with $\Ga$).

\begin{remark}\rm
Replacing the containment $X\in \GL(m,\bbZ)$ by $X\in
\GL(m,\bbR)$, we come to a criterion for $M$ and $\td M$ to be
almost casually isometric. Indeed, the manifolds are almost
isometric if and only if they admit realizations with equal
parameters, except for $\Gamma$, but all lattices in $\bbR^m$ are
linearly equivalent. For $m=1$ and $a\neq 0$, we have an ordinary
quadratic polynomial,which is equivalent to $A+s^2$, $A>0$. If
$k=0$, then $\dim M=4$, and the main theorem implies that
nonelliptic manifolds in this dimension form a one-parameter
family. There is a more precise version of this result
\cite[Theorem~1]{GM}: all flat complete strictly casual
nonelliplic Lorentzian 4-manifolds are nomothetic to the manifold
of the following example.
\end{remark}
\begin{example}\rm
Let $\Ga=\bbZ$ be the cyclic infinite group generated by the
affine transformation $v\to\la v+\tau$ which is defined in the
basis $e_0,\dots,e_3$ by relations
\begin{eqnarray*}\label{fourm}
&\la e_0=e_0,\quad\la e_1=e_1,\quad\la e_2=e_2+e_0,\quad\la
e_3=e_3+e_2+\frac12 e_0,\quad\tau=e_1;\quad{}\\
&\ell(v,v)=2v_0v_3-v_1^2-v_2^2.\nonumber
\end{eqnarray*}
The manifold $M=V/\Gamma$ is strictly casual. Let
$u=(u_0,\dots,u_3)\in V$, $u_3>0$. A straightforward calculation
(see \cite{GM}) shows that the past of $u$ contains the open
halfplane $x_3<-\frac1{u_3}$, meets the hyperplane
$x_3=-\frac1{u_3}$, but does not include it. In particular, this
implies that the past of $u$ contains a straight line and is not
closed.
\end{example}

Our aim is now to describe the quadratic polynomials $Q(s)$ in
(\ref{poabc}) such that $Q=Q_M$ for some $M$. We will achieve in
the two steps: in the first, we reduce the problem to the case of
nondegenerate matrix $C$; in the second, we describe the
polynomials with $C>0$.
\begin{proposition}\label{degen}
Let $(n,m,r,k)$ be the signature of a manifold $M$ and let $k>0$.
Then there exists $X\in\GL(m,\bbR)$ such that the matrix $X^\top
Q_M(s)X$ admits block-diagonal realization with blocks of size
$k\times k$ and $(m-k)\times(m-k)$, where $k$-block does not
depend on $s$ and $(m-k)$-block is a characteristic polynomial for
a manifold $\td M$ of the signature $(n-k,m-k,r,0)$; $M$ is almost
causally isometric to the product of $\td M$ and a flat $k$-torus.
Moreover,
\begin{eqnarray}\label{rankc}
m-k=\rank C.
\end{eqnarray}
\end{proposition}
\begin{theorem}\label{descc}
The polynomial $Q(s)$ in (\ref{poabc}), where $A,B,C$ are
$m$-matrices, defines a characteristic curve of manifold $M$ with
signature $(n,m,r,0)$ if and only if  (\ref{poabc}) holds,
$m+r+2\leq n$, and
\begin{eqnarray}\label{cccri}
&C-BA^{-1}B\geq0,\\
&\label{detr} r=\rank (C-BA^{-1}B).
\end{eqnarray}
If $m=r$ then (\ref{poabc}) may be replaced by a weaker condition
$A>0$.
\end{theorem}
It follows from (\ref{poabc}), (\ref{cccri}), and (\ref{detr})
that $r>0$ (see remark at the end of the paper). If (\ref{cccri})
is true and $A>0$, then (\ref{poabc}) can be formulated in term of
the eigenvectors of some matrices as in (\ref{freec}).

\begin{remark}\rm
 The conditions
(\ref{cccri}) is not a consequence of (\ref{poabc}) as well as
(\ref{poabc}) does not follow from (\ref{cccri}) even under the
additional assumption $A>0$. For example, let $m=2$, $A=B=\one$,
and
\begin{eqnarray*}
C=\left(
\begin{array}{cc}
1&\varepsilon\\
\varepsilon&1
\end{array}
\right),\qquad 0<\varepsilon<1.
\end{eqnarray*}
Then (\ref{poabc}) is true but (\ref{cccri}) is false. If $m=2$,
$A=C=\one$, $B$ is diagonal with entries 1 and 0 in the diagonal
then $Q(-1)$ is degenerate; hence (\ref{poabc}) is false but
(\ref{cccri}) is true in this case.
\end{remark}
The condition (\ref{cccri}) means that $Q(s)$ is, roughly
speaking, a sum of squares (see the proof of the theorem). If
$m=1$, then $BA^{-1}B-C$ is the discriminant (divided by $A$) of
the quadratic polynomial $Q(s)$.

The characteristic curve implicitly contains some invariants of
$M$. An essential instance is given in the following proposition.
It does not determine $M$ completely, in particular, it does not
distinguish almost isometric and homothetic manifolds. The
spectrum of a matrix $X$ is denoted by $\spp(X)$.
\begin{proposition}\label{spectr}
Let polynomials $Q_M(s)=A+2sB+s^2C$ and $Q_{\td M}(s)=\td A+2s\td
B+s^2\td C$ be characteristic for almost causally isometric
manifolds $M$, $\td M$, and let $C>0$, $\td C>0$. Then eigenvalues
of $BC^{-1}$, $\td B\td C^{-1}$ are real and
\begin{eqnarray}\label{invsp}
\spp\big(BC^{-1}\big)=\al\big(\spp\big(\td B\td C^{-1}\big)\big)
\end{eqnarray}
for some $\al\in\Aff(\bbR)$.
\end{proposition}
\begin{remark}\rm
Using results above, it is not difficult to describe the flat
complete strictly causal manifolds in small dimensions. For $n=4$
there is exactly one, up to a homothety, nonelliptic manifold (see
(see Example~1). Let $n=5$, $e_1,\dots,e_{5}$ be the standard
basis for $V$,
\begin{eqnarray*}
\ell(u,u)=2u_{4}u_{5}-u_1^2-u_2^2-u_{3}^2.
\end{eqnarray*}
If $m+r+2=5$ and $k=0$, then $m=2$ and $r=1$. By (\ref{freec}),
$a'$ must have a simple spectrum. Putting $ae_1=t e_3$,
$ae_2=e_2+r e_3$, where $t,r\neq0$, $v_0=e_5$, $v_1=e_4$, and
choosing $\Ga$ in the linear span of $e_1,e_2$, we get all
manifolds of the signature $(5,2,1,0)$. Other signatures can be
reduced to less dimensions. The number of variants grows rather
fast with $n$.
\end{remark}

\section{Proof of results}
\begin{proof}[Proof of Lemma~1] Since $ax\not \in L$,
$\lambda(x)v=v$ and (\ref{forla}) imply $\l_0(v)=0$. Hence, $v\in
W$. Further, we have $W\cap(C\cup(-C))=L$ since $W$ is tangent to
$\partial C$ at $v_0$. The same relations imply the second
assertion of (1) in the lemma. Since $v_0,x,ax\in W$, all
hyperplanes $W_s$ are $\Gamma$-invariant by (\ref{forla}) and
(\ref{forta}). Putting $v_0=0$ and $l_0(v)=s$ in (\ref{forla}) and
(\ref{forta}), we obtain the formula for the action of $\Gamma$ in
$W_s/L$: $\gamma_x(v)=v+sax+x (\mod L)$. Thus, $T$ acts by pure
translations in $W_s/L$, and the action is free and proper if
(\ref{freea}) is true. Then it is free and proper in $W_s$. This
proves (3) and (2) of the lemma.
\end{proof}
\begin{proof}[Proof of Proposition~\ref{zarcl}]
Let $p$ be a polynomial on $V$.  Suppose that $p(\ga_x(v))$ is
independent of $x\in\Ga$ for some $v\in V$. It follows from
(\ref{forla}) and (\ref{forta}) that $p(\ga_x(v))$ is a polynomial
on $x$. Hence, $p(\gamma_x(v))$ is constant on $T$. Therefore, the
Zariski closure of $\Ga v$ includes $Tv$. To prove the reverse
inclusion, note that the projection of a $T$-orbit to $V/L$ is the
affine subspace $\phi((\one+sa)T+v)$, where $s=l_0(v)$. Hence,
$T$-orbit has codimension 1 in the affine subspace
$X=L+(\one+sa)T+v$. By Lemma~\ref{trivi}, (3), it has the form
$\{v+x+sax+f(x)v_0:\,x\in T\}$ for some function $f$ on $T$. One
can find $f$ by simple straightforward calculation with
(\ref{forla}) and (\ref{forta}), which shows also that the orbit
is distinguished by an algebraic equation in $X$. Thus, every
$T$-orbit is Zariski-closed. Clearly, the image of $T$ under the
embedding defined by (\ref{forla}) and (\ref{forta}) is also
closed. Since the action of $T$ is free by (\ref{freea}), the same
is true for $\wh \Gamma$.
\end{proof}
\begin{proof}[Proof of Theorem~\ref{uniho}]
It follows from (\ref{forla})
that
\begin{eqnarray*}
(\la(x)-1)V\subset W, \quad(\la(x)-1)W\subset L,
\quad(\la(x)-1)L=0.
\end{eqnarray*}
Hence, every $M(a,\Ga)$ has  unipotent holonomy group. To prove
the converse, note that $M$ can be finitely covered by $M(a,\Ga)$
for some $a,\Ga$ by \cite[Theorem~2]{GM}. Then $M=V/\td\Ga$, where
$\td\Ga\subset\cP(V)$ is an unipotent group that contains $\Ga$ as
a subgroup of finite index and acts in $V$ freely and properly.
Then $\wh\Ga$ has finite index in the algebraic closure
$\wh{\td\Ga}$ of $\td\Ga$. Each unipotent matrix $U$ lies in a
one-parameter group $\exp(tX)$, where $X$ is nilpotent and is a
polynomial of $U-\one$; moreover, $\exp(tX)$ is polynomial in $t$.
Hence, the Zariski closure of the cyclic group generated by $U$
includes $\exp(\bbR X)$. Therefore, every Zariski-closed unipotent
linear group is connected. This implies that $\wh{\td\Ga}=\wh\Ga$
and $\td\Ga\subset\wh\Ga$. By Proposition~\ref{zarcl}, $\td\Ga$ is
a discrete subgroup of the vector group $\wh\Ga$ including the
uniform lattice $\Ga$. Hence, $\td\Ga$ is a uniform lattice
itself; thus, $M=M(a,\td\Ga)$.
\end{proof}
\begin{proof}[Proof of Proposition~\ref{geome}]
By Theorem~\ref{uniho}, we may assume $M=M(v_0,v_1,T,a,\Ga)$.
Suppose first that $p=\ka(o)$. Then we may identify $V$ and
$T_pM$. It follows from (\ref{forta}) that $d_0\tau_(x)=x$.
Hence,the tangent space at $p$ to the orbit of $T$ coincides with
$T$ under this identification. By Proposition~\ref{zarcl}, the two
definitions of $T$ agree. Since $M$ is nonelliptic, $ax\neq0$ for
some $x\in T$. Furthermore, $\ell(ax,ax)\neq0$ because $\ell$ is
negative definite on $T$. By (\ref{forla}), for all sufficiently
large $t>0$ and all $v\in V$ such that $l_0(v)\neq 0$, we have
$(\la(tx)v-v)+(\la(-tx)v-v)=rv_0$, where $r\neq0$. Hence, $v_0\in
H$ and we see that $H=aT+\bbR v_0$ as an immediate consequence of
(\ref{forla}). Since $v_0$ is isotropic and $\ell$ is
nondegenerate on $aT$, $L=H\cap H^\bot=\bbR v_0$ and $\dim L=1$.
Let  $\td v_0=\frac1r v_0$ for some $r>0$, $\td v_1= r v_1+w$,
where $w\in T^\bot\cap W$ is such that $\td v_1$ is isotropic.
Then, according to (\ref{forla}),
\begin{eqnarray*}
\la(x)\td v_1=\la(x)(r v_1+w)=rax+\xi(x,w,r)v_0
\end{eqnarray*}
for some function $\xi$. Put $\td a=ra$. Then
$\ell(ax,x)v_0=\ell(\td ax,x)\td v_0$. Since the kernel of the
orthogonal projection in $V$ to every subspace complementary to
$L$ in $W$ always intersects $W$ by $L$ and $ax\in U$, we come to
the formulas (\ref{lwn})--(\ref{forta}) for the action of $T$,
with $v_0,v_1$, and $a$ replaced by $\td v_0,\td v_1$, and $\td
a$. The embedding of $\Ga$ to $T$ satisfies the equality
\begin{eqnarray*}
\Ga=\pi_T(\exp_p^{-1}(p))
\end{eqnarray*}
and is completely determined by the latter. Evidently,
$\exp_p(\ga(v))=\exp_p(v)$ for all $v\in V$ and $\ga\in\Ga$.
Further, the inclusion $E\subseteq N$ implies that $\ell$ is
nondegenerate in $E$. Hence, $V=E\oplus E^\bot$. The decomposition
is $\Ga$-invariant since $E^\bot\supseteq T+aT+L$. Therefore, $E$
is a direct factor in $M=V/\Ga$. Since the action of $\Gamma$ in
$E^\bot$ is subject to the same formulas as in $V$, it follows
that $M'$ satisfies the proposition.

To prove the assertion for any $p\in M$, it is sufficient to
remove origin to an arbitrary point $\td o\in V$ and find
parameters that realize the action in the form
(\ref{lwn})--(\ref{forta}). A pure translation does not change the
linear parts of affine transformations. Therefore, $H$ and $L$ are
the same as above. By the first part of the proof, we may preserve
$v_1$. Then $N$ does not change. Consequently,
\begin{eqnarray*}
&\la(x)=\td\la(\td x),\\
&\td a\td x =\pi_{N}\td\la(\td x)v_1=\pi_{N}\la(x)v_1=ax
\end{eqnarray*}
(the tilde distinguishes new parameters), where $\td x$ is the
point in the tangent space $\td T$ to the orbit of $\td o$ at $\td
o$ that satisfies  $\phi(\td x)=\phi(\td\tau(\td x))$, where
\begin{eqnarray}\label{tdtau}
\td\tau(\td x)=\ga_x(\td o)-\td o=\tau(x)+(\la(x)-\one)\td o,\quad
x\in T.
\end{eqnarray}
Differentiating by $x$ at $x=0$ we find
\begin{eqnarray*}
&\td T=\{x+sax-\ell(ax,\td o)v_0:\,x\in T\},\quad s=l_0(\td o);\\
&\td x=x+sax-\ell(ax,\td o)v_0.
\end{eqnarray*}
The latter is true since $\phi(\td x)=\phi(x+sax)$ by
(\ref{tdtau}) and (\ref{forla}), (\ref{forta}). If $s=0$
(equivalently, if $\td o\in W$), then $\phi(\td T)=\phi(T)$;
$U=T+H$ does not change since $L\subset H$. Inserting this into
(\ref{forla}), (\ref{forta}), we obtain the same formulas with new
parameters. Let $\td o=sv_1$, where $s\in\bbR\setminus\{0\}$. Then
$l_0(\td o)=s$, $\ell(ax,\td o)=0$, and we find
\begin{eqnarray*}
&\td\tau(\td x)=\tau(x)+(\la(x)-\one)\td
o=x-\frac12\ell(ax,x)v_0+sax -\frac{s}2\ell(ax,ax)v_0=\\
&\td x-\frac12\ell(\td a\td x,\td x)v_0.
\end{eqnarray*}
Since each translation in $V$ is a composition of a translation
along $W$ and $\bbR v_1$, the action can be realized in the form
(\ref{lwn})--(\ref{forta}) for any choice of origin.
\end{proof}

\begin{proof}[Proof of Theorem~\ref{main}]
The assertion on signatures is clear. The left-hand side of
(\ref{defqv}) is the squared length of the unique segment in $V$
that represents $x\in\pi_1(M,p)$. Hence, an isometry identify
$q_v$, $\td q_v$ as functions on $\pi_1(M,p)$. Identifying $\Ga$
with $\bbZ^n$, we arrive at quadratic polynomials with matrix
coefficients by (\ref{raspi}) and (\ref{qsabc}). Clearly, the
transformations of the coefficients which are induced by a change
of parameters are subject to (\ref{eqcur}). If images of
polynomials $Q_M$ and $Q_{\td M}$ (which can be a parabola, a ray,
or a single point) coincide, then there exists an increasing
affine change of the variable $s$ that identifies them. Since any
change of this type can be realized by shifting the origin along
$\bbR v_1$ and scaling the vectors $v_0$, $v_1$ (see the proof of
Proposition~\ref{geome} above and (\ref{shifs}), (\ref{affch})),
we arrive at equal curves in $P_m/\GL(m,\bbZ)$ if the manifolds
are causally isometric. Hence, (\ref{eqcur}) is true.

Conversely, let (\ref{eqcur}) hold. It is sufficient to prove the
existence of a transformation in $\cP_n$ equivariant with respect
to the action of the fundamental groups of $M$ and $\td M$ in $V$
assuming that they are subject to (\ref{lwn})--(\ref{forta}).
Applying a transformation in $\cP(V)$, we may assume $v_0=\td
v_0$, $v_1=\td v_1$ (hence, $N=\td N$ and $W=\td W$); and also
that $T=\td T$ and $R=\td R$ (since the signatures coincide).
Thus, it is sufficient to prove that $a$ and $\td a$ are
conjugated by an $\ell$-orthogonal linear transformation of $N$
preserving $T$ and $R$. Using (\ref{raspi}) and the decomposition
(\ref{deca}), by (\ref{defqv}) and (\ref{qsabc}) we find
\begin{eqnarray}
&\scal{Az}{z}=\ell(x,x),\label{ax}\\
&\scal{Bz}{z}=\ell(a'x,x),\label{bx}\\
&\scal{Cz}{z}=\ell(a'x,a'x)+\ell(a''x,a''x),\label{cx}
\end{eqnarray}
where $x=\iota\,z\in T$, $z\in\bbR^m$; similar equalities hold for
$Q_{\td M}$. Thus, for instance, (\ref{ax}) means that
$\ell(x,x)=\ell(\td x,\td x)$, where $\td x=\td\iota z$. It
follows from (\ref{ax}) that $\iota =\xi\td\iota$ for some
$\xi\in\OO(\ell|_{T})$. Then $\td a'=\xi^{-1}a'\xi$ by (\ref{bx})
(note that the quadratic form on the left-hand side uniquely
determines $a'$ and $\td a'$ since they are $\ell$-symmetric).
This implies $\ell(a'x,a'x)=\ell(\td a'\td x,\td a'\td x)$ for all
$x\in T$. Then $\ell(a''x,a''x)=\ell(\td a''\td x,\td a''\td x)$
by (\ref{cx}). Therefore, $\td a''=\ze^{-1} a''\xi$ for some
$\ze\in\OO(\ell|_{R})$. Thus, the transformation which is equal to
$\xi$ on $T$, equal to $\ze$ on $R$, and is identical on $U^\bot$,
identifies parameters for the two actions (including the embedding
of $\Ga$).
\end{proof}
\begin{proof}[Proof of Proposition~\ref{degen}]
Since $a$ is $\ell$-symmetric, $\ker a\perp aT$. The form $\ell$
is negative definite on $N$; hence $\ker a\cap aT=0$. There exists
$X\in\GL(m,\bbR)$ such that $\iota X\iota^{-1}$ identifies the
decomposition $T=\ker a\oplus\td T$, where $\td T=T\cap(\ker
a)^\bot$, with the natural decomposition:
$\bbR^m=\bbR^k\oplus\bbR^{m-k}$. It follows from (\ref{raspi})
that $q_s(x)$ is independent of $s$ if $x\in\ker a$. Put $\td
V=(\ker a)^\bot$ and $\td\Ga=\bbZ^{m-k}$. Since $\td V$ contains
$\td T$, $aT=a\td T$, and $L$; it is $\td\Ga$-invariant by
(\ref{lwn})--(\ref{forta}), and the action is subject to the same
formulas. Put $\td M=\td V/\td\Ga$. Clearly, $M$ is almost
causally isometric to the product of $\td M$ and the torus
$\bbR^k/\bbZ^k$ (note that $\ga_x(v)=v+x$ if $x\in \ker a$).
Comparing the coefficients of $s^2$ in (\ref{raspi}) and
(\ref{qsabc}), we find $\rank C=\rank a$. This proves
(\ref{rankc}). Remaining assertions are obvious.
\end{proof}
In what follows, the fractional powers of nonnegative matrices are
supposed nonnegative.
\begin{lemma}\label{spech}
Let $Q$ be as in (\ref{poabc}), $C>0$. Any transformation
(\ref{linch}) with $X\in\GL(m,\bbR)$ does not change
$\spp\big(C^{-\frac12}BC^{-\frac12}\big)$ and
$\spp\big(A^{-\frac12}BA^{-\frac12}\big)$.
\end{lemma}
\begin{proof}
If two quadratic polynomials of the type $A+2sB+s^2\one$ are
conjugated by a transformation $X$ as in (\ref{linch}), then $X$
is orthogonal. Therefore,  the spectrum of the coefficient of $s$
does not depend on the choice of the polynomial with coefficient
$\one$ at $s^2$ in the $\GL(m,\bbR)$-orbit under the
transformations (\ref{linch}) of the polynomial $A+2sB+s^2C$.
Analogous arguments show that spectra of the coefficients of $s$
of all polynomials of the type $\one+2sB+s^2C$ in this orbit are
equal. Thus, the equalities
\begin{eqnarray}\label{prive}
&C^{-\frac12}(A+2sB+s^2C)C^{-\frac12}=
C^{-\frac12}AC^{-\frac12}+2sC^{-\frac12}BC^{-\frac12}+s^2\one,\\
&A^{-\frac12}(A+2sB+s^2C)A^{-\frac12}=
\one+2sA^{-\frac12}BA^{-\frac12}+s^2A^{-\frac12}CA^{-\frac12}\nonumber
\end{eqnarray}
prove the lemma.
\end{proof}
\begin{proof}[Proof of Proposition~\ref{spectr}] Spectra of
$C^{-\frac12}BC^{-\frac12}$ and $BC^{-1}$ are identical since
these matrices are conjugated; they are real because the first
matrix is symmetric. Replacing $s$ by $t(s-s_0)$ in $A+2sB+s^2C$,
we find coefficients $t^2C$ and $2(tB-t^2s_0C)$ at $s^2$ and $s$,
respectively. So, $BC^{-1}$ corresponds to
$\frac1tBC^{-1}-s_0\one$. It remains to apply Theorem~\ref{main}
and Lemma~\ref{spech}.
\end{proof}
\begin{proof}[Proof of Theorem~\ref{descc}]
Let $Q(s)=Q_M(s)$ for some $M$. Put $\td
B=A^{-\frac12}BA^{-\frac12}$, $\td C=A^{-\frac12}CA^{-\frac12}$.
Then
\begin{eqnarray}\label{squa}
Q(s)=A^{\frac12}(\one+2s\td B+s^2\td
C)A^{\frac12}=A^{\frac12}\left((\one+s\td B)^2+s^2(\td C-\td
B^2)\right)A^{\frac12}.
\end{eqnarray}
On the other hand,
$q_s(x)=\ell((\one+sa')x,(\one+sa')x)+s^2\ell(a''x,a''x)$.
Comparing this with (\ref{squa}), we find
$-\ell(a''x,a''x)=\scal{(\td C-\td B^2)z}{z}$, where
$x=\iota\,A^{\frac12}z$, $z\in\bbR^m$.
Therefore,
\begin{eqnarray}
&C-BA^{-1}B=A^{\frac12}\left(\td C-\td
B^2\right)A^{\frac12}\geq0,\label{nonne}\\
&r=\rank a''=\rank(\td C-\td B^2)=\rank(C-BA^{-1}B).\nonumber
\end{eqnarray}
The condition (\ref{freea}) holds if and only if
$\rank(\one+sa)=m$ for all $s\in\bbR$; in other words, it is true
if and only if the form on the left-hand side of (\ref{raspi}) is
positive definite, which is equivalent to (\ref{poabc}).

Conversely, let $Q(s)$ satisfy (\ref{poabc}) and (\ref{cccri}).
Then the manifold $M$ may be constructed following (not word for
word) the procedure (A)--(C)  Proposition~\ref{geome}.
\begin{itemize}
\item[(A)] Put $V=\bbR^{n}$, where $n\geq m+r+2$ and $r$ is
defined by (\ref{detr}). Let $e_1,\dots,e_n$ be the standard basis
for $\bbR^n$,
$$\ell(z,z)=2z_{n}z_{n-1}-z_1^2-\dots-z_{n-2}^2,$$
$v_0=e_{n-1}$, $v_1=e_{n}$, and let $L,W,N$ be as in (\ref{lwn}).
Define special subspaces that were introduced above by the
following decomposition:
\begin{eqnarray*}
\bbR^n=\bbR^m\oplus\bbR^r\oplus\bbR^{n-m-r-2}\oplus\bbR\oplus\bbR=
T\oplus R\oplus E\oplus\bbR v_1\oplus L.
\end{eqnarray*}
\item[(B)] Put $a'=\td B$ and $a''=J(\td C-\td B^2)^{\frac12}$,
where $J$ is a linear isometry of the range of $\td C-\td B^2$
onto $R$ . \item[(C)] Define $\iota:\,\bbR^m\to T$ as
$\iota=A^{-\frac12}$.
\end{itemize}

The condition $m=r$ is equivalent to $C-BA^{-1}B>0$; if $A>0$ then
$\td B$ and $\td C$ are well defined. Furthermore, (\ref{nonne})
taken together with (\ref{squa}), implies (\ref{poabc}).
\end{proof}
It follows from (\ref{squa}) that the equality $C-BA^{-1}B$ hold
if and only if $Q(s)$ admits representation in the form $X(\one
+sY)^2X$, where $X$ and $Y$ are symmetric matrices. If $Y\neq 0$
then $Q(s)$ is degenerate for some $s\in \bbR$. Thus,
(\ref{poabc}), (\ref{cccri}),and (\ref{detr}) imply that $r>0$.

\end{document}